\documentstyle[12pt,psfig]{article}

\input{psfig.sty}

\newtheorem{theorem}{Theorem}
\newtheorem{lemma}[theorem]{Lemma}
\newtheorem{proposition}[theorem]{Proposition}

\newtheorem{corollary}[theorem]{Corollary}

\newtheorem{definition}[theorem]{Definition}

\newtheorem{problem}[theorem]{Problem}

\newenvironment{proof}{{\it Proof:\/}}{$\Box$\vskip 0.08in}

\begin{document}

\begin{center}
\begin{LARGE}
\baselineskip=10pt {\bf Non-left-orderable 3-manifold groups
}
\end{LARGE}
\end{center}
\begin{center} Mieczys{\l}aw K.~D{\c a}bkowski,
J\'ozef H.~Przytycki and Amir A.~Togha
\end{center}

\vspace{27pt}
\centerline{Abstract}
{\footnotesize{ We show that several torsion free 3-manifold groups 
are not left-orderable. 
Our examples are groups of cyclic branched coverings of $S^3$
branched along links.
 The figure eight knot provides simple 
nontrivial examples. The groups arising in these examples are known
as Fibonacci groups which we show not to be left-orderable. 
Many other examples of non-orderable groups are obtained by taking 
3-fold branched covers of $S^3$ branched along various hyperbolic 
2-bridge knots. 
The manifold obtained in such a way from the $5_2$ knot
is of special interest as it is conjectured to be the hyperbolic 
3-manifold with the smallest volume.
}}
\ \\
\ \\
We investigate the orderability properties of fundamental groups 
of 3-dimensional manifolds. 
We show that several torsion free 3-manifold groups
are not left-orderable. Many of our manifolds are obtained by taking
n-fold branched covers along various hyperbolic 2-bridge knots.
The paper is organized in the following way:\ after defining 
left-orderability we state our main theorem listing branched set links
and multiplicity of coverings from which we obtain manifolds with 
non-left-orderable groups. Then we describe presentations of these groups
in a way which allows the proof of non-left-orderability in a uniform way.
The Main Lemma (Lemma 5) is the algebraic underpinning of our method 
and the non-left-orderability follows easily from it in almost all cases.
Moreover we prove  the non-left-orderability of a family of 
3-manifold groups to
which the Main Lemma does not apply. These groups, known as
generalized Fibonacci groups $F(n-1,n)$, arise as groups of double 
covers of $S^3$ branched along pretzel links of type $(2,2,...,2,-1)$.
We end the paper with some  questions and speculations.

\begin{definition}\label{1}
A group is {\textup {left-orderable}} if there is a strict total ordering 
$\prec$ of its elements which is left-invariant: $x \prec y$ iff $zx \prec zy$
for all $x,y$ and $z$.
\end{definition}

Straight from the definition, it follows that a group with a torsion 
element is not left-orderable.

It is known that groups of compact, $P^2$-irreducible 3-manifolds with
non-trivial first Betti number are left-orderable \cite{BRW,H-S}.
However, our main theorem below lists various classes of 
 3-manifolds with non-left-orderable groups. Non-left-orderability 
of 3-manifold groups has interesting consequences for the geometry
of the corresponding manifolds \cite{C-D,RSS}.

\begin{theorem}\label{2}
Let $M_L^{(n)}$ denote the $n$-fold  branched cover of $S^3$
branched along the link $L$, where $n >1$. Then the fundamental group,
$\pi_1(M_L^{(n)})$, is not left-orderable in the following cases:
\begin{enumerate}
\item[(a)] $L=T_{(2',2k)}$ is the torus link of the type $(2,2k)$ with the
anti-parallel orientation of strings, and $n$ is arbitrary (Fig.1).

\item[(b)] $L=P(n_1,n_2,...,n_k)$ is the pretzel link of the type 
$(n_1,n_2,...,n_k)$, $k>2$, 
where either (i) $n_1,n_2,...,n_{k} > 0$, or 
(ii) $n_1=n_2= \cdots =n_{k-1}=2$ and $n_k = -1$ (Fig.2). 
The multiplicity of the covering is $n=2$.

\item[(c)]$L=L_{[2k,2m]}$ is a 2-bridge knot of the type $\frac{p}{q}=
2m + \frac{1}{2k}=[2k,2m]$, where $k, m >0$, and $n$ is arbitrary (Fig.4). 
\item[(d)] $L=L_{[n_1,1,n_3]}$ is the 2-bridge knot of the type $\frac{p}{q}=
n_3 + \frac{1}{1 + \frac{1}{n_1}} = [n_1,1,n_3]$, 
where  $n_1$ and $n_3$ are odd positive numbers. 
The multiplicity of the covering is $n\leq 3$.
\end{enumerate}
\end{theorem}
\ \\
\centerline{\psfig{figure=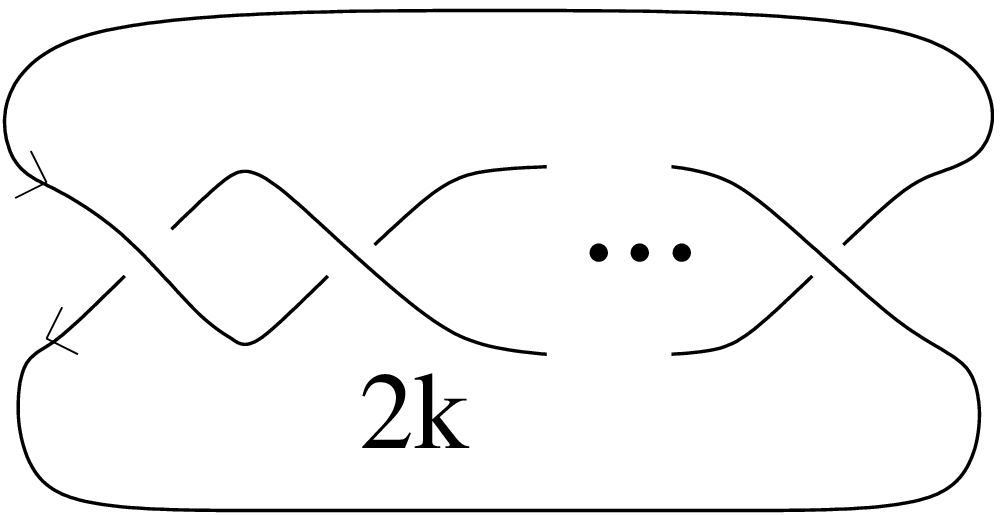,height=1.9cm}}
\begin{center} Fig. 1 \end{center}
\ \\
\ \\
\centerline{\psfig{figure=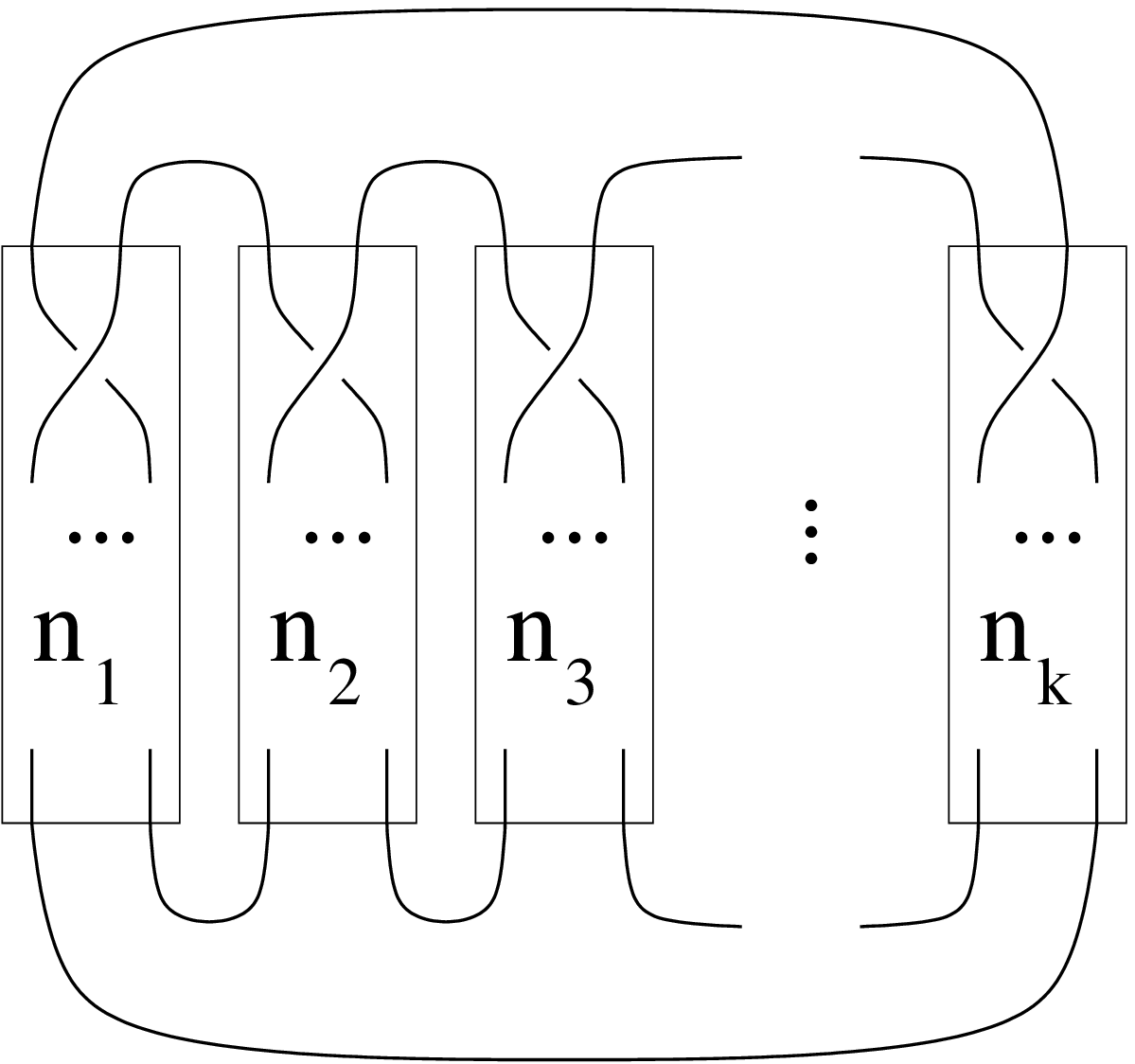,height=5.1cm}}
\begin{center} Fig. 2 \end{center}

The manifolds described in parts (a), (b), and also for $n\leq 3$ and 
the figure eight knot, $L=L_{[2,2]}=4_1$, in part (c) 
are Seifert fibered manifolds. The non-left-orderability of
their groups follows from the general characterization of Seifert fibered 
manifolds with a left-ordering \cite{BRW}. 
Part (c) for the figure eight knot when $n=3$ is of historical 
interest because it was the first known example of 
a non-left-orderable torsion free 3-manifold 
group \cite{Rol}\footnote{This Euclidean 
manifold was first considered by
Hantzsche and Wendt \cite{H-W}. J.~Conway has proposed to call this
manifold {\it didicosm}. It can be also described as the 2-fold branched cover
over $S^3$ branched along the Borromean rings.}. 
Part (c) for the figure eight knot when $n>3$, 
gives rise to hyperbolic manifolds that are
related to examples discussed in \cite{RSS}, as they are Dehn fillings 
of punctured-torus bundles over $S^1$. The manifolds obtained in parts
(c) and (d), when $n>2$ (except $M^{(3)}_{4_1}$), are all hyperbolic 
manifolds as well\footnote{It follows from the Orbifold 
Theorem that branched $n$-fold covers ($n>2$) of $S^3$ branched along 
hyperbolic 2-bridge knots and links or along the Borromean rings are 
hyperbolic, except for $M^{(3)}_{4_1}$ which is a Euclidean manifold, 
didicosm \cite{Bo,HLM,Ho,Th}.}. 

The case $\frac{p}{q}= \frac{7}{4} =1 + \frac{1}{1+\frac{1}{3}}=[3,1,1]$, 
that is, the branching set being the $5_2$ knot, 
is of special interest since $M^{(3)}_{5_2}$ is conjectured 
to be the hyperbolic 3-manifold with the smallest volume \cite{Ki}. 
The fact that $\pi_1(M^{(3)}_{5_2})$
is not left-orderable was observed in \cite{C-D,RSS}. 
The non-left-orderability in other cases is proved here for the first time.

The special form of the presentations of the groups listed in Theorem 2,
allows us to conclude the theorem in most cases, using the 
Main Lemma formulated below (Lemma 5).
\begin{proposition}\label{3} 
The groups listed in Theorem 2 have the following presentations:\\
\textup{(}a\textup{)} 
$\pi_{1}(M^{(n)}_{T_{(2',2k)}}) =\\  \{x_{1},x_{2},...,x_{n} |\ \
x^{k}_{1}x^{-k}_{2}=e,\ x^{k}_{2}x^{-k}_{3}=e,...,
   x^{k}_{n}x^{-k}_{1}=e, \ x_{1}x_{2}\cdots x_{n}=e\}$ \\
\textup{(}b\textup{)} 
\textup{(}i\textup{)} $\pi_{1}(M^{(2)}_{P_{(n_{1},n_{2},...,n_{k})}}) =\\ 
\{x_{1},x_{2},...,x_{k} | \ \ x_{1}^{n_1}x_{2}^{-n_2}=e, 
x_{2}^{n_2}x_{3}^{-n_3}=e,
...,x_{k}^{n_k}x_{1}^{-n_1}=e ,\ x_1x_2\cdots x_k=e\}$\\
\ \ \ \ \textup{(}ii\textup{)} $\pi_{1}(M^{(2)}_{P_{(2,2,...,2,-1)}}) = \{
x_{1},x_{2},...,x_{k} | \ \ x_{1}^{2}=x_{2}^{2}= \cdots =x^{2}_{k}=x_{1}x_{2}
\cdots x_{k}\}$ \\ 
\vspace*{0.1cm}
\textup{(}c\textup{)}
$\pi_{1}(M^{(n)}_{L_{[2k,2m]}}) =\\
\vspace*{0.1cm}
 \{z_{1},z_2,\ldots,z_{2n} | \ \ 
z_{2i+1}= z_{2i}^{-k}z_{2i+2}^k,\ 
z_{2i} = z_{2i-1}^{-m}z_{2i+1}^m,\ 
z_2z_4\ldots z_{2n}=e \}$  where $i=1,2 \ldots n$  and   
subscripts are taken modulo $2n$. \\
\textup{(}d\textup{)} 
$\pi_{1}(M^{(n)}_{L_{[2k+1,1,2l+1]}}) = \{x_{1},...,x_{n} | \ \
r_{1}=e, ..., r_{n}=e,\  
x_{1}x_{2}\cdots x_{n}=e\}$, where $k \geq 0$, $l \geq 0$, $$r_{i} =
x_{i}^{-1}(x_{i}^{-k}x_{i+1}^{k+1}x_{i}^{-1})^{l}x_{i}^{-k}x_{i+1}^{k+1}
((x_{i+1}^{
-k}x_{i+2}^{k+1}x_{i+1}^{-1})^{l}x_{i+1}^{-k}x_{i+2}^{k+1})^{-1},$$ 
and subscripts are taken modulo $n$.
\end{proposition}

\begin{proof}
Since the presentations for all manifolds from Theorem 2 are obtained by
similar calculations, therefore we shall only provide
full details for the case (c) (compare \cite{M-V}). 
Let $T_{1}$ denote the 2-tangle in Fig.3(a),  $-[2k]$ in Conway's notation 
and let $T_{2}$ denote the 2-tangle in Fig.3(b),  $[2m]$ in Conway's notation.
Let us assume that the arcs of $T_{1}$ and $T_{2}$ 
are oriented in the way shown in Fig.3.\\
\ \\
\centerline{\psfig{figure=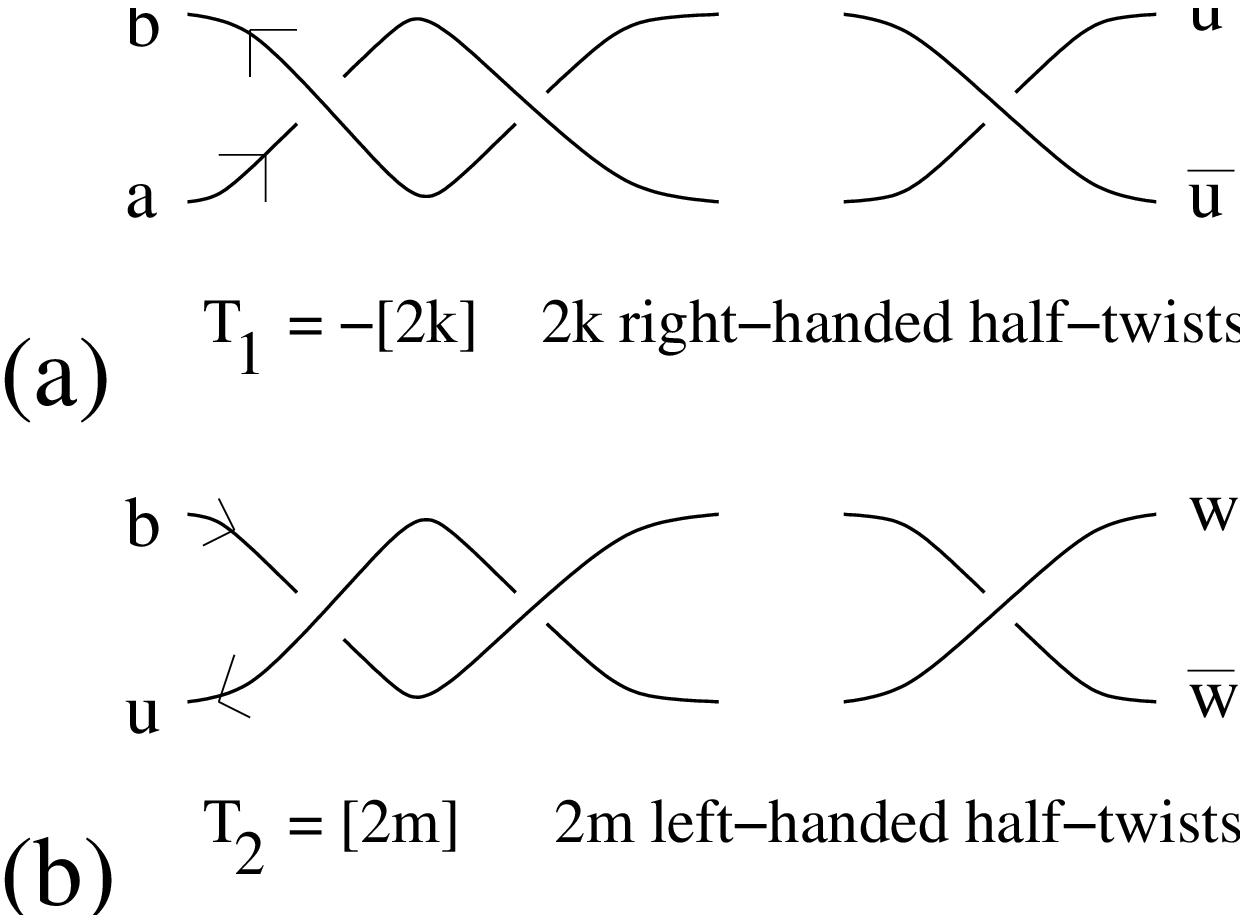,height=4.4cm}}
\begin{center} Fig. 3 \end{center}

Let $F_{2} = \{a,b \ | \  \- \}$ be a free group generated by $a$ and $b$.
Assign to the initial arcs of $T_1$ the generators $a$ and $b$.
Then by successive use of Wirtinger relations, progressing from 
left to right in the diagram, we finally decorate the terminal arcs  
by $\bar{u} = (ba^{-1})^{k}a(ab^{-1})^{k}$ and
$u = (ba^{-1})^{k}b(ab^{-1})^{k}$, 
respectively (see Fig.3(a)).
Analogously, assigning to initial arcs of the tangle
$T_{2} = [2m]$ (Fig.3(b)) the elements $b$ and $u$ of $F_{2}$ and using
Wirtinger relations successively one obtains terminal arcs decorated by
$w = (u^{-1}b)^{m}b(b^{-1}u)^{m}$ and $\bar{w} = 
(u^{-1}b)^{m}u(b^{-1}u)^{m}$, respectively.
\ Combining these calculations in the fashion illustrated in Fig.4, we 
obtain the relation $((ba^{-1})^{k}b^{-1}(ab^{-1})^{k}b)^{m}b= 
a((ba^{-1})^{k}b^{-1}(ab^{-1})^{k}b)^m$ and the presentation
$$\pi_{1}(S^{3}-L_{[2k,2m]}) = \{a,b | \ \
r=((ba^{-1})^{k}b^{-1}(ab^{-1})^{k}b)^{m}b
((ba^{-1})^{k}b^{-1}(ab^{-1})^{k}b)^{-m}a^{-1} \}.$$ 

\ \\
\centerline{\psfig{figure=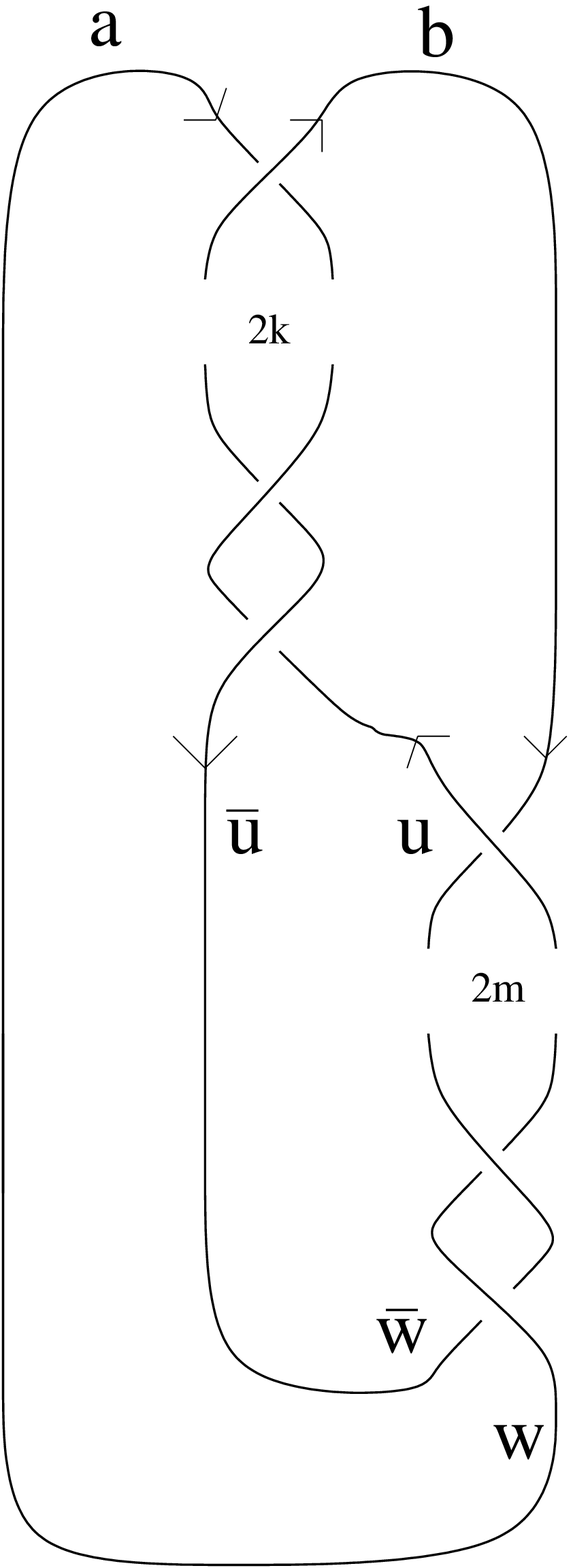,height=8.7cm}}
\begin{center} Fig. 4; The 2-bridge knot $[2k,2m]$ \end{center}

Using Fox non-commutative calculus \cite{Cr}, as explained in
 \cite{Pr,P-R}, we compute a presentation of $\pi_{1}(M^{(n)}_{L_{[2k,2m]}})$ 
by ``lifting" the generators $a$ and $b$ as well as the defining relation 
$r$ of $\pi_{1}(S^{3}-L_{[2k,2m]})$.

We illustrate this by first computing a presentation of the fundamental 
group of the $n$-fold cyclic {\it unbranched} covering of $S^3-L_{[2k,2m]}$. 
Since $\pi_{1}(S^{3}-L_{[2k,2m]})$ has 2 generators, $a$ and $b$, 
the covering space  
will have $n+1$ generators, that is, $y=ab^{-1}, \tau(y),\tau^2(y)...,
\tau^{n-1}(y)$ and $b^n$, where $\tau$ is the inner automorphism of $F_2$, 
given by $\tau(w)=bwb^{-1}$ (see Fig.5). \\
\begin{figure}[h]
\centerline{\psfig{figure=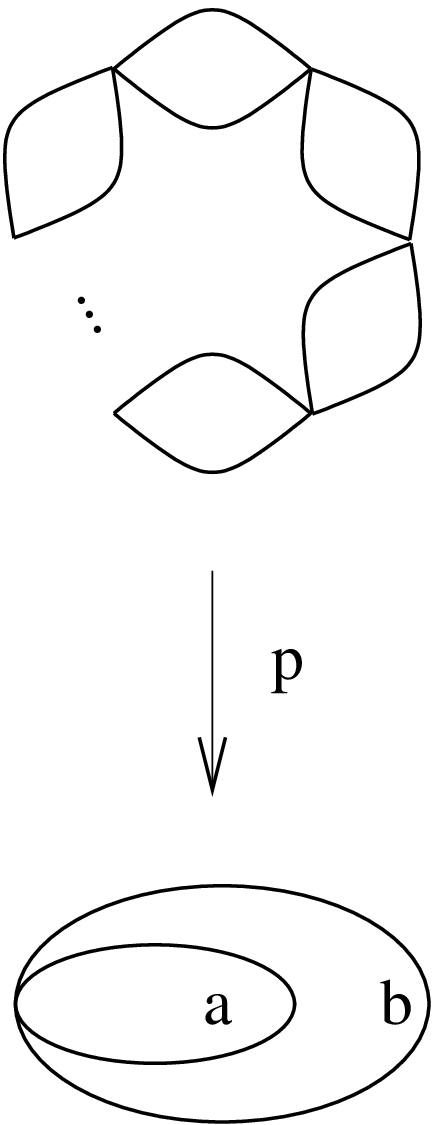,height=5.1cm}}
\centerline{Fig. 5}
\end{figure}

The relation $r$ will also be lifted to $n$ relations $\tilde r, 
\tau(\tilde r), \tau^2(\tilde r),...,\tau^{n-1}(\tilde r)$, in 
the group of the $n$-fold cyclic covering, where 
$$\tilde r = 
((y^{-1})^{k}(\tau^{-1}(y))^{k})^{m}((\tau(y^{-1}))^{k}(y)^{k})^{-m} y^{-1}.$$

When dealing with the branched case, however, the relations 
$a^n=e$ and $b^n=e$ should also be added\footnote{Since $L_{[2k,2m]}$ 
is a knot, 
 the relation $a^n=e$ follows from the relation $b^n=e$ and the relations 
$\tau^i(\tilde r)$.}.
We then write the word $a^n$ in terms of new generators as 
$y\tau(y)...\tau^{n-1}(y)$. 
In order to simplify the presentation of 
$\pi_{1}(M^{(n)}_{L_{[2k,2m]}})$ we put 
$x_{1}=\tau^{-1}(y),\, 
x_{2}=y,\,x_{3}=\tau(y),\,\dots,\, x_{n}=\tau^{n-2}(y)$.
Thus 
$$\pi_{1}(M^{(n)}_{L_{[2k,2m]}}) = \{x_{1},\,x_{2},\dots,\, x_{n} | \ \
x_{i}^{-1}(x_{i}^{-k}x_{i-1}^{k})^{m}(x_{i+1}^{-k}x_{i}^{k})^{-m} = e
,\ x_{1}x_{2}\cdots x_{n}=e \},$$
where $i = 1,\,2,\,\dots,\,n$ and subscripts are taken modulo $n$.

To change this presentation to the one described in Proposition 3(c)
we ``deform" variables by putting
$z_{2i}= x_i$ and $z_{2i+1}= x_{i}^{-k}x_{i+1}^k$. In new variables 
the presentation has the desired form 
$$\pi_{1}(M^{(n)}_{L_{[2k,2m]}}) =
\{z_{1},z_2,\ldots,z_{2n} | \ \
z_{2i+1}= z_{2i}^{-k}z_{2i+2}^k,
z_{2i} = z_{2i-1}^{-m}z_{2i+1}^m,
z_2z_4\cdots z_{2n}= e \},$$ 
where $i=1,2,...,n$ and subscripts are taken modulo $2n$.
\footnote{In the special case of $k=m=1$ we obtain the classical 
Fibonacci group $F(2,2n)$ already known to be the fundamental group of 
$M^{(n)}_{4_1}$. We suggest that the presentation 
for any $k$ and $m$ to be called the ($k,m$)-{\it deformation}, 
$F((k,m),2n)$, of the classical Fibonacci group.}
\end{proof}
It is worth mentioning that the case (c) that we singled out for illustrating
the proof of Proposition 3 involves a step that the proofs for other cases do
not require. More specifically, all of the presentations given in the
statement of Proposition 3, except for the case (c), are results of
straightforward calculations and we do not need to deform
the variables in any way in those cases in order 
to obtain the desired presentation.


The following definition and Main Lemma capture the algebraic 
properties of listed groups.
\begin{definition}\label{4}
\begin{enumerate}
\item[(i)]
Given a finite sequence $\epsilon_{1}, \epsilon_{2}, ..., \epsilon_{n}$,
$\epsilon_{i} \in \{ -1, 1\}$, for all $i = 1, 2, ..., n$
and a nonempty reduced word $w =
x_{a_{1}}^{b_{1}}x_{a_{2}}^{b_{2}}...x_{a_{m}}^{b_{m}}$ of the
free group
$F_{n} = \{ x_{1}, x_{2}, ..., x_{n} \ | \  \}$, 
we say $w$ \textup {blocks} the sequence 
$\epsilon_{1}, \epsilon_{2}, ..., \epsilon_{n}$
if either $\epsilon_{a_{j}}b_{j} > 0$ for all $j$ or
$\epsilon_{a_{j}}b_{j} < 0$ for all $j=1, 2, ..., m$. 
\item[(ii)]
A set $W$ of reduced words of $F_{n}$ is {\textup {complete}} 
if for any given
sequence $\epsilon_{1}, \epsilon_{2}, ..., \epsilon_{n}$,
$\epsilon_{i} \in \{ -1, 1\}$, for $i = 1, 2, ..., n$, there is a word 
$w \in W$ that blocks $\epsilon_{1}, \epsilon_{2}, ..., \epsilon_{n}$.
\item[(iii)] The presentation $\{ x_{1}, x_{2}, ..., x_{n}\ \ | \  W\}$ of
a group $G$ is called {\textup {complete}} if the set $W$ of relations
is complete.
\end{enumerate}
\end{definition}

\begin{lemma}[Main Lemma]\label{5}
Any nontrivial group $G$ that admits a complete presentation is
not left-orderable.
\end{lemma}

\begin{proof}
Suppose, on the contrary, that $\prec$ is a left-ordering on $G$.  
Let $G = \{ x_{1}, x_{2}, ..., x_{n}\ \ | \  W\}$ be a complete 
presentation of $G$. 
Let 
$E = \{ (\epsilon_{1}, \epsilon_{2}, ..., \epsilon_{n})\  | \ 
x_{i}^{\epsilon_{i}} \preceq e$\ in the group $G$, where 
$\epsilon_{i} \in \{ -1, 1\}, i = 1, 2, ..., n \}$. 
Since $W$ is complete, each sequence $(\epsilon_{1}, \epsilon_{2}, ...,
\epsilon_{n}) \in E$ is blocked by a word $w \in W$.
Since $w$ is a relator, this is impossible, because the product 
of a number of ``positive" elements in a left-orderable group 
will be ``positive", not the identity.
This contradiction completes the proof.
\end{proof}

Theorem 2 follows easily from the Main Lemma  and Proposition 3 
in all cases except
for part (b)(ii) which we deal with separately in the following lemma.
\begin{lemma}\label{6}
Let $F(n-1,n) = $ \ \\
$\{x_1,\cdots,x_n$ $\vert$
$x_{1}x_{2} \cdots x_{n-1}=x_n,\ 
x_{2}x_{3} \cdots x_{n}=x_1,\ 
\cdots, x_{n}x_{1} \cdots x_{n-2}=x_{n-1}\}$. \\
If $n>2$, then $F(n-1,n)$ is not left-orderable.
\end{lemma}
\begin{proof}
 $F(2,3)$ is finite (it is the quaternion group $Q_8$),
hence it is not left-orderable.
Let us assume, then, that $n>3$. First of all, note that the mapping 
$x_i \mapsto g:
F(n-1,n) \rightarrow \{g$ $\vert$ $g^{n-2}=e\}=Z_{n-2}$ 
defines an epimorphism, and since $n-2>1$
our group is not the trivial group. 

It is not hard to see that in $F(n-1,n)$ we have 
$x^{2}_1=x_{2}^{2}= \cdots =x^{2}_n=
x_{1}x_{2} \cdots x_{n}$. 
Let $t=x^{2}_i=x_{1}x_{2} \cdots x_{n}$ for any $i$.
Suppose
that $\prec$ is a left-ordering on $F(n-1,n)$. Since
$F(n-1,n)$ is not the trivial group, hence $t \neq e$ unless 
our group has a torsion, which is not the case.
Consider the case $t \prec e$. The case $e \prec t$ can be dealt with 
similarly.

Since $t=x_i^2$, 
we must have $x_i \prec e$ for all $i$. In particular, $x_i \neq
e$ for all $i$. 
This makes $x_1 \preceq x_2 \leq \cdots \preceq x_n \preceq x_1$ impossible,
because if
$x_1=x_2= \cdots =x_n \neq e$, then  
$x_1^2=t=x_{1}x_{2} \cdots x_{n}=x_1^n$
implies
$x^{n-2}_1=e$, which in turn makes $F(n-1,n)$ a torsion group and thus
non-left-orderable.

Therefore, $x_{i+1} \prec x_{i}$ for some $i$ modulo $n$. 
Assume, without loss of generality,
that $x_n \prec x_{n-1}$. Multiplying from the left
by
$x_{1}x_{2} \cdots x_{n-1}$ one obtains
$$t=x_{1}x_{2} \cdots x_{n-1}x_n \prec x_{1}x_{2} \cdots x_{n-2}x_{n-1}x_{n-1} =
x_{1}x_{2} \cdots x_{n-2}t = tx_{1}x_{2} \cdots x_{n-2} .$$
The last equality holds because $t=x_i^2$ commutes will all $x_i$. Multiplying
both
sides from the left by $t^{-1}$ gives $e \prec x_{1}x_{2} \cdots x_{n-2}$,
contradicting
the fact that $x_{i} \prec e$ for all $i$. 
\end{proof}

Left-orderability of a countable group $G$ is equivalent to $G$ being
isomorphic to a subgroup of ${Homeo}_+({\bf R})$ (compare \cite{BRW}). 
Calegari and Dunfield related left-orderability of the group 
of a 3-manifold $M$ with foliations on $M$. Therefore 
we have
\begin{corollary}\label{7}
\begin{enumerate}
\item[(i)] The groups of manifolds described in Theorem 2
do not admit a faithful representation to $Homeo_+({\bf R})$.
\item[(ii)] Manifolds described in Theorem 2 do not admit
a co-orientable ${\bf R}$-covered foliation {\textup{\cite{C-D}}}.
\end{enumerate}
\end{corollary}
Thurston proved that if
an atoroidal 3-manifold $M$ has a taut foliation 
then there exists a faithful  action of $\pi_1(M)$ 
on $S^1$ \cite{C-D}. Exploring the fact that the group of the 
manifold of the smallest known volume,
$M_{5_2}^{(3)}$, (together with some of its subgroups) is not
left-orderable, Calegari and Dunfield showed that $\pi_1(M_{5_2}^{(3)})$
does not admit a faithful  action of $\pi_1(M)$  on $S^1$ and therefore
$M_{5_2}^{(3)}$ does  not admit a taut foliation \cite{C-D}.
The connection between  faithful  actions of $\pi_1(M)$  on $S^1$ 
and on ${\bf R}$ is to be explored further. 

  We would like to contrast our non-left-orderability results with 
some examples of left-orderable 3-manifold groups.


It is known that if $M_{K}^{(n)}$ is irreducible 
(as is always the case for a hyperbolic knot $K$) 
and the group $H_1(M_{K}^{(n)})$ is infinite,  
then the group 
$\pi_1(M_{K}^{(n)})$ is left-orderable \cite{BRW,H-S}. 
 There are several examples of 2-bridge knots with infinite
homology groups of cyclic branched coverings along them.
For the trefoil knot $3_1$ we have $H_1(M_{3_1}^{(6k)}) = Z\oplus Z$. 
For hyperbolic 2-bridge knots $9_6=K_{[2,2,5]}$ and
$10_{21} = K_{[3,4,1,2]}$ the groups $H_1(M_{9_6}^{(6)})$ and 
$H_1(M_{10_{21}}^{(10)})$ are also infinite \footnote{To see quickly that 
$H_1(M_{K}^{(n)})$ is infinite one can use Fox theorem which says that 
$H_1(M_{K}^{(n)})$ is infinite if and only if the Alexander polynomial,
 $\Delta_K(t)$, is equal to zero for some $n$th root 
of unity. To test the last condition 
for small knots one can use tables of knots with $\Delta_K(t)$ 
decomposed into irreducible factors \cite{B-Z}. We check, for example, 
that ${\Delta}_K(e^{\pi i/3})=0$ for hyperbolic 2-bridge knots
$K= 8_{11}, 9_6, 9_{23}, 10_5, 10_9, 10_{32}$ and $10_{40}$. 
Note also that Casson and Gordon proved that $p^k$-fold cyclic 
branched coverings
along a knot, where $p$ is prime, are rational homology spheres.}.

We end the paper with some questions about possible
generalizations of our results.
\begin{problem}\label{8}
\begin{enumerate}
\item[(i)] Are the groups  $\pi_1(M_{5_2}^{(n)})$ non-left-orderable
for $n>3$?
\item[(ii)] Are the groups  $\pi_1(M_{K}^{(n)})$ of
hyperbolic 2-bridge knots $K$  with finite  $H_1(M_{K}^{(n)})$
 non-left-orderable?
\item[(iii)] Are the groups  $\pi_1(M_{K}^{(n)})$ of
hyperbolic knots $K$ with finite $H_1(M_{K}^{(n)})$
 non-left-orderable?
\item[(iv)] In general, for which links $L$ and multiplicities of
covering $n$, is the group $\pi_1(M_L^{(n)})$ non-left-orderable?
\end{enumerate}
\end{problem}
\ \\
{\bf Acknowledgment.}  We would like to thank Andrzej Szczepa\'nski for
introducing Fibonacci groups to us. We are also grateful to 
 Jos\'e Montesinos, Dan Silver and 
Andrei Vesnin for their valuable correspondence.

\ \\
Department of Mathematics\\
The George Washington University \\
Washington, DC 20052 \\
e-mails: \\
\texttt{mdab@utdallas.edu} \\
\texttt{przytyck@gwu.edu}\\
\texttt{userid@gwu.edu}
\end{document}